\documentclass[11pt,twocolumn,twoside]{article}
\usepackage{fully3d}

\usepackage{float}
\usepackage{textcomp}
\usepackage{tikz}
\usepackage[export]{adjustbox}
\usetikzlibrary{spy}
\usepackage{tabularx}
\usepackage{amssymb,amsmath,amsfonts, amsthm}
\usepackage{caption}
\usepackage{subcaption}

\usepackage{mathtools}
\usepackage[percent]{overpic}
\usepackage{algorithm,algcompatible,amsmath}
\algnewcommand\INPUT{\item[\textbf{Input:}]}%
\algnewcommand\PARAMETER{\item[\textbf{Parameters:}]}%
\algnewcommand\OUTPUT{\item[\textbf{Output:}]}%

\newcommand{\R}{\mathbb R}

\newcommand{\XX}{\mathbf{x}} 
\newcommand{\PP}{\mathbf{p}}
\newcommand{\QQ}{\mathbf{q}}
                      
\newcommand{\ZZ}{\mathbf{z}}

\newcommand{\Ad}{\mathbf A}

\newcommand{\trans}{{\scriptstyle \boldsymbol{\mathsf{T}}}}

\newcommand{\LLambda}{\boldsymbol{\Lambda}}

\begin{document}

%-------------------------------------------------------------------------------------------
%%%%% add your title here %%%%%
%-------------------------------------------------------------------------------------------
%%%%% add your title here %%%%%
\title{ Unrolled three-operator splitting for parameter-map learning \\ in Low Dose X-ray CT reconstruction } 

%%%%% add authors and affiliations here %%%%%
\author[1]{Andreas Kofler}
\author[2,3]{Fabian Altekrüger}
\author[3]{Fatima Antarou Ba}
\author[1]{Christoph Kolbitsch}
\author[4,*]{Evangelos Papoutsellis}
\author[1]{David Schote}
\author[5]{Clemens Sirotenko}
\author[1]{Felix Frederik Zimmermann}
\author[6]{Kostas Papafitsoros}

\affil[1]{Physikalisch-Technische Bundesanstalt (PTB), Braunschweig and Berlin, Germany~~~~~~~~~~}
\affil[2]{Humboldt-Universit\"at zu Berlin, Department of Mathematics, Berlin, Germany~~~~~~~~~~~~~~~}
\affil[3]{Technische Universit\"at Berlin, Institute of Mathematics, Berlin, Germany~~~~~~~~~~~~~~~~~~~~~}
\affil[4]{Finden Ltd, Rutherford Appleton Laboratory,  Harwell Campus, Didcot, United Kingdom~~~~}
\affil[5]{Weierstrass Institute for Applied Analysis and Stochastics, Berlin, Germany~~~~~~~~~~~~~~~~~~}
\affil[6]{School of Mathematical Sciences, Queen Mary University of London, United Kingdom~~~~~\vspace{1em}}
\affil[*]{Corresponding author: epapoutsellis@gmail.com}

%%%%% don't change these 2 lines %%%%%
\maketitle
\thispagestyle{fancy}

%-------------------------------------------------------------------------------------------
%%%%% add your summary (abstract) here               %%%%%%
%%%%% use footnotesize for this section              %%%%%%
%%%%% please stick to the customabstract environment %%%%%% 

\begin{customabstract}
We propose a method for fast and automatic estimation of spatially dependent regularization maps for total variation-based (TV) tomography reconstruction. 
%Our work is inspired by recent developments in algorithm unrolling using deep neural networks (NNs), the Primal Dual Three-Operator Splitting (PD3O) is used for the total variation (TV) tomography reconstruction problem. 
The estimation is based on two distinct sub-networks, with the first sub-network estimating the regularization parameter-map from the  input data while the second one unrolling $T$ iterations of the Primal-Dual Three-Operator Splitting (PD3O)  algorithm. The latter approximately solves the corresponding TV-minimization problem incorporating the previously estimated regularization parameter-map. The overall network is then  trained end-to-end in a supervised learning fashion using pairs of clean-corrupted data but crucially without the need of having access to labels for the optimal regularization parameter-maps.     
\end{customabstract}

%-------------------------------------------------------------------------------------------
%%%%% main text                                                %%%%%    
%%%%% remove the dummy content and put your own content here   %%%%% 
%%%%% feel free to choose your own section titles              %%%%% 
%%%%% you don't need to put the content in a separate tex file %%%%%

% dummy_content.tex shows how to add sections, figures, tables, formulas, and references
% remove the following line, it just adds dummy content
\section{Introduction}

Over recent years, Low Dose X-ray Computed Tomography (LDCT) has received a growing interest in the medical imaging field due to its ability to reduce the radiation dose. Patients are exposed to low levels of radiation by reducing the energy of the photons emitted from the X-ray source. Using traditional and analytic reconstruction methods such as filtered back projection (FBP), several imaging artifacts are introduced, compromising the quality of the reconstructed image and clinical diagnosis.

To overcome this problem, iterative reconstruction methods have been proposed such as algebraic reconstruction technique (ART), simultaneous algebraic reconstruction technique (SART) and projection onto convex sets (POCS). In addition, such reconstruction procedures often require the use of regularization methods in order to eliminate noise and artifacts, such as for instance, the well-known Tikhonov and Total Variation (TV) regularization \cite{sidky2012convex}.

The acquired measured tomography data can be described by the  equation $ \ZZ = \Ad \XX_{\mathrm{true}} + \mathbf{e}$,
where $\XX_{\mathrm{true}}\in \mathbb{R}^{n}$ is the ground truth image, $\Ad: \mathbb{R}^{n}\to \mathbb{R}^{m}$ is a linear operator which models the data-acquisition process, i.e. the discretized Radon transform, and  $\mathbf{e}\in \mathbb{R}^{m}$ denotes some random noise component. 
% The goal is to reconstruct $\XX_{\mathrm{true}}$ or at least a good enough approximation of it given the data $\ZZ$. 
Regularized iterative methods solve minimization problems of the form
\begin{equation}\label{eq:reg_problem}
    \underset{\XX }{\min}  \, \mathcal{D}(\Ad \XX, \ZZ) + \mathcal{R}(\XX),
\end{equation}
where $\mathcal{D}(\,\cdot \,, \, \cdot \,)$ denotes a data-discrepancy measure and $\mathcal{R}(\,\cdot\,)$ a regularization term. A classical example is the TV tomography reconstruction problem under Gaussian noise which can be written as 
\begin{equation}\label{eq:tv_min_problem}
    \underset{\XX }{\min}\; \frac{1}{2} \| \Ad \XX - \ZZ\|_2^2 + \lambda \| \nabla \XX \|_1 + \mathbb{I}_{\{\XX>0\}}(\XX).
\end{equation}
%It is well known that although TV regularization eliminates noisy artifacts, it also introduces staircasing artifacts, e.g., 
%blocky-like, piecewise constant structures. This can be improved using high order methods such as the total generalized variation (TGV), \cite{TGV}. However, 
A key factor which impacts the quality of the reconstructed image is the careful choice of the regularization parameter $\lambda$ which balances the strength between the regularization and the data fidelity term. Underestimating $\lambda$ yields poor regularization, while overestimating it results in smooth images with an artificial ``cartoon-like'' appearance. Particularly in medical imaging applications, where images are at the basis of diagnostic decisions and therapy planning, a proper choice of any regularization parameter is crucial. 

Employing a single scalar parameter $\lambda$ implies that the regularization is enforced with equal strength for each pixel/voxel. Depending on the application, this might be undesirable due to different features contained in the image. In this case, one can replace the scalar parameter with a spatially varying, i.e.\ a pixel/voxel dependent one, denoted now by $\boldsymbol{\Lambda}\in \R_{+}^{qn}$. Here, $q$ denotes the number of directions for which the partial derivatives are computed. Implementation-wise, $\boldsymbol{\Lambda}$ corresponds to a stack of diagonal operators which contain a different regularization parameter for every single pixel/voxel in the respective gradient domain of the image. Then, the resulting problem  has the form
\begin{equation}\label{eq:tv_min_problem_ldct}
\underset{\XX }{\min}\; \mathcal{D}(\Ad \XX , \ZZ) +\|\boldsymbol{\Lambda} \nabla \XX\|_{1} + \mathbb{I}_{\{\XX>0\}}(\XX).
\end{equation}
However, this problem requires a precise data-adaptive estimation of the spatially varying parameter-map $\boldsymbol{\Lambda}$  which is a highly non-trivial task.

\begin{figure*}[h!]
    \centering
    \includegraphics[width=0.8\textwidth]{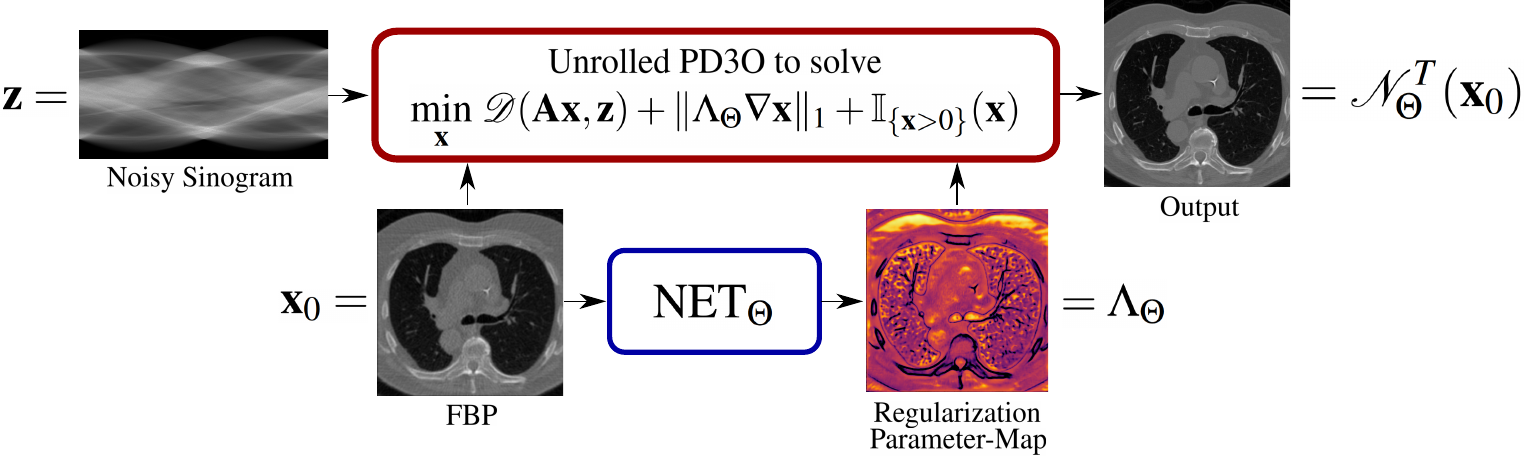}
    \caption{\small{ Network architecture for the LDCT reconstruction problem. It consists of a sub-network $\mathrm{NET}_{\Theta}$ that estimates the regularization parameter-map (blue), and a sub-network that unrolls $T$ iterations of  the PD3O algorithm (red).  
    }}
    \label{fig:tv_pd3O_nn.pdf}
\end{figure*}

\section{Methods}
One approach  for the automatic estimation of the spatially varying regularization parameter $\LLambda$ is employing  bilevel optimization techniques. Given $M$ pairs of measured data  and the corresponding ground truth $(\ZZ_{i},\XX_{\mathrm{true}}^{i})_{i=1}^{M}$, the general bilevel formulation is
\begin{equation}\label{general_bilevel}
\left \{
\begin{aligned}
&\min_{\boldsymbol{\Lambda}}\; \sum_{i=1}^{M}l(\XX^{i}(\boldsymbol{\Lambda}),\XX_{\mathrm{true}}^{i}) \;\text{subject to }\;\;\\ & \XX^{i}(\boldsymbol{\Lambda})=\underset{\XX}{\operatorname{argmin}}\; \mathcal{D}(\Ad \XX , \ZZ) + \| \boldsymbol{\Lambda}\nabla \XX \|_1  + \mathbb{I}_{\{\XX>0\}}(\XX),
\end{aligned}\right.
\end{equation}
where $l$ is a suitable upper level objective.  For instance, if $l(x_{1},x_{2})=\|x_{1}-x_{2}\|_{2}^{2}$, the bilevel problem \eqref{general_bilevel} seeks to compute the parameters $\boldsymbol{\Lambda}$ which are ``the best on average'', i.e. PSNR-maximizing, for the given $M$ data pairs. Hence, given some new data $\ZZ_{\mathrm{test}}$ which has been measured in a similar way as $(\ZZ_{i})_{i=1}^{M}$, solving \eqref{eq:tv_min_problem_ldct} with the precomputed $\boldsymbol{\Lambda}$ will yield a good reconstruction.

 Although this scheme has been extensively studied both for scalar and spatially varying regularization parameters, it has been mainly applied in image denoising applications, i.e. $\Ad=\mathbf{I}_n$ with Gaussian noise, see \cite{pockbilevel, bilevellearning, Chung_Reyes_Schoenlieb}.  Further, unsupervised approaches employing upper level energies that do not depend on the ground truth $\XX_{\mathrm{true}}$, i.e., $l:=l(\XX(\boldsymbol{\Lambda}))$, have also been  considered in a series of works \cite{hintermuellerPartII,  bilevel_handbook, bilevelTGV,  bilevelconvex}. Even though these bilevel optimization methods are typically supported by rigorous mathematical theories, they are computationally demanding which has limited their use on tomographic problems.

\subsection{An Unrolled Neural Network Framework}

Here, inspired by the recent success of unrolled neural networks (NNs)  \cite{Monga_2021}, we consider an unrolled neural network 
approach in order to learn the regularization parameter $\boldsymbol{\Lambda}$. The proposed framework  is   summarized in Figure \ref{fig:tv_pd3O_nn.pdf} and it is outlined next.

%When a scalar parameter is employed, the authors in \cite{Afkham_2021}, follow a supervised learning approach for computerized tomography reconstruction and image deblurring. Similarly to the bilevel optimization framework, the pipeline consists of an offline and an online phase. In the first phase a family of optimal regularization parameters $(\lambda_{i})_{i=1}^{M}$ is computed, e.g.\ by employing a scheme like \eqref{general_bilevel}. Then, in the second part of the offline phase,  using the training data $D=\{(\lambda_{i}, \ZZ_{i})_{i=1}^{M}\}$, the parameters $\Theta$  of a NN $\mathcal{N}_{\Theta}$ are learned by minimizing
%\begin{equation}\label{sup_learn_param_off1}
%\min_{\Theta} \mathcal{L}(\Theta) := \frac{1}{M} \sum_{i=1}^{M} l( \mathcal{N}_{\Theta}(\ZZ_{i}),\lambda_{i}),
%\end{equation}
%for a suitable loss function $l$. Once an estimate of the optimal parameters $\Theta$ has been learned, one passes to the online phase, and given some new data $\ZZ_{\mathrm{test}}$, the regularization parameter is simply calculated by applying the learned network to $\ZZ_{\mathrm{test}}$, i.e., $\lambda_{\Theta}=\mathcal{N}_{\Theta}(\ZZ_{\mathrm{test}})$ \eqref{eq:tv_min_problem} is solved by an appropriate algorithm.

%Inspired by the recent success of unrolled NNs  \cite{Monga_2021}, a different strategy is applied for the construction of the regularization parameter-maps. 
An unrolled NN which corresponds to an implementation of an iterative scheme  of finite length is constructed to approach the solution of problem \eqref{eq:tv_min_problem} assuming a \textit{fixed} regularization parameter-map. Within the unrolled NN, the regularization parameter-map is estimated from the input data via a sub-network $\mathrm{NET}_{\Theta}$ and is used throughout the whole reconstruction scheme.  To be more precise, given some initial estimate $\XX_{0}$ we work with an iterative scheme (speficied in the next section)
\begin{equation}\label{iteraive_scheme}
\XX_{T}=S^{T}(\XX_{0},\ZZ, \boldsymbol{\Lambda},\Ad), \quad T=0,1,2,\ldots,
\end{equation}
that solves \eqref{eq:tv_min_problem_ldct} in the limit $T\to\infty$.
Then, for some fixed number of iterations $T\in \mathbb{N}$, our unrolled NN reads as follows:
\begin{equation}\label{unrolled_intro}
\left \{
\begin{aligned}
\boldsymbol{\Lambda}_{\Theta}&=\mathrm{NET}_{\Theta}(\XX_0),\\
\XX_{1}&=S^{1}(\XX_{0},\ZZ, \boldsymbol{\Lambda}_{\Theta},\Ad),\\
\vdots\\
\XX_{T}&=S^{T}(\XX_{0},\ZZ, \boldsymbol{\Lambda}_{\Theta},\Ad).
\end{aligned}
\right.
\end{equation}
Here, $\mathrm{NET}_{\Theta}$ denotes a U-Net \cite{Ronneberger2015} with learnable parameters $\Theta$. We denote by $\mathcal{N}_{\Theta}^{T}$ the overall resulting network, i.e. 
\[\mathcal{N}_{\Theta}^{T}(\XX_0)=S^{T}(\XX_{0},\ZZ, \boldsymbol{\Lambda}_{\Theta},\Ad)
=S^{T}(\XX_{0},\ZZ, \mathrm{NET}_{\Theta}(\XX_0),\Ad).\] The unrolled NN can then be end-to-end trained in a supervised manner on a set of input-target image-pairs. This resulting network can be identified as a pipeline that combines in a sequential way 1) the estimation of the regularization parameter-map which is adapted to the  data $\ZZ$ (and hence in medical imaging to the new patient) and 2) the iterative scheme that solves the  image reconstruction problem.

\subsection{Primal-Dual Three-Operator Splitting}

The iterative scheme selected here for the LDCT reconstruction problem is the Primal-Dual Three-Operator (PD3O) splitting algorithm. The PD3O was introduced in \cite{Yan2018} and it is a generalized version of the Primal-Dual Hybrid Gradient (PDHG) algorithm \cite{chambolle2011first}.  It is used to minimize  objectives that consist of a proximable function $g$, a composite function $f$ with the linear operator $\mathbf{K}$ and a differentiable function $h$ with a Lipschitz constant $L$:
\begin{align*}
    \min_\XX f(\mathbf{K}\XX) + g(\XX) + h(\XX).
\end{align*}
\begin{algorithm}[t]
\caption{Unrolled PD3O algorithm}\label{algo:general_pd3o_algo}
  \begin{algorithmic}[1]
  \INPUT $L = \text{Lip}(\nabla h)$,\; $\tau = 2/L$, \; $\sigma = 1/(\tau \Vert \mathbf{K} \mathbf{K}^\trans\| )$, \; $\text{initial guess}~ \bar{\XX}_0$
  \OUTPUT reconstructed image $\XX_{\mathrm{TV}}$
  \STATE $\PP_0 = \mathbf{\bar{x}_0}$
  \STATE $\QQ_0 = \mathbf{0}$
    \FOR {$k < T$ }
    \STATE $\QQ_{k+1} = \mathrm{prox}_{\sigma f^{*}}(\QQ_k + \sigma \mathbf{K}\bar{x}_k)$
    \STATE $\PP_{k+1} = \mathrm{prox}_{\tau g}(\PP_k - \tau \nabla h(\PP_k) - \tau \mathbf{K}^\trans \QQ_{k+1})$
    \STATE $\bar{\XX}_{k+1} = 2\PP_{k+1} - \PP_k + \tau\nabla h(\PP_k) - \tau \nabla h(\PP_{k+1})$ 
    \ENDFOR
    \STATE $\XX_{\mathrm{TV}} = \XX_T$
  \end{algorithmic}
\end{algorithm}

The algorithm is summarized in Algorithm \ref{algo:general_pd3o_algo} and explained next. Unlike the standard $\mathrm{L}^{2}$-squared fidelity term that is commonly used in tomography reconstruction problems with Gaussian noise, here we employ the Kullback-Leibler divergence which is more suitable to describe the noise distribution of the measured tomographic data $\ZZ$. We have that
$\ZZ= \Ad \XX + \mathbf{e}$, where
\begin{align*}
 \mathbf{e} = -\Ad \XX - \log (\mathbf{\tilde{N}_1} / N_0) ~,\, \mathbf{\tilde{N}_1} \sim \text{Pois}(N_0 \exp(-\Ad \XX \mu)).
\end{align*}
We denote with $\mu$ and $N_0$  the normalization constant and the mean photon count per detector bin without attenuation, respectively. The data-discrepancy in \eqref{eq:tv_min_problem_ldct} can be derived from a Bayesian viewpoint and is 
\begin{align} \label{eq_CTdatafidelity}
\mathcal{D}(\Ad \XX ,\ZZ) = \sum_{i=1}^m e^{-(\Ad \XX )_i \mu} N_0 - e^{-\ZZ_i \mu} N_0 \big(-(\Ad \XX)_i \mu + \log(N_0) \big),
\end{align}
see \cite{ADHHMS2022} for more details. To configure PD3O for \eqref{eq:tv_min_problem_ldct} we define the following
\begin{align*}
& f(\QQ) = \|\boldsymbol{\Lambda} \QQ\|_{1}, \quad g(\PP) = \mathbb{I}_{\{\PP>0\}}(\PP), \quad \mathbf{K} = \nabla,\\
& h(\PP) = \sum_{i=1}^m e^{-\PP_i \mu} N_0 - e^{-\ZZ_i \mu} N_0 \big(-\PP_i \mu + \log(N_0) \big).
\end{align*}

Notice that with the standard $\mathrm{L}^{2}$-squared fidelity term, it is sufficient to use the  PDHG algorithm since its convex conjugate has a closed-form proximal operator, which is not the case with  \eqref{eq_CTdatafidelity}. However, the additional function in the PD3O algorithm allows to  express the data discrepancy in the differentiable term $h$. Note that $\nabla h$ is not globally Lipschitz continuous but due to the non-negativity constraint, we only have to consider $\nabla h(\PP)$ for $\PP$ with non-negative entries. Consequently, we can find an upper bound of the Lipschitz constant of $\nabla h$ by $\text{Lip}(\nabla h) \le \Vert \Ad \Vert^2 \mu^2 N_0$.

\section{Results}

To evaluate our proposed unrolled NN, we use the LoDoPaB dataset \cite{LoDoPaB21} for low-dose CT imaging. It is based on scans of the Lung Image Database Consortium and Image Database Resource Initiative which serve as ground truth images, while the measurements are simulated. The dataset contains 35820 training images, 3522 validation images and 3553 test images. Here the ground truth images have a resolution of $362\times 362$ on a  domain of $26\text{cm}\times26\text{cm}$. We only use the first 300 training images and the first 10 validation images. For the forward operator we consider a normalization constant $\mu = 81.35858$, the mean photon count per detector bin $N_0 = 4096$ as well as 513 equidistant detector bins and 1000 equidistant angles between 0 and $\pi$.

\begin{figure}[t!]
\centering
\begin{subfigure}[t]{.5\textwidth}
\centering
\caption*{\textbf{FBP}}
  \begin{tikzpicture}[spy using outlines={rectangle,white,magnification=2,size=1.5cm, connect spies}]
\node[anchor=south west,inner sep=0]  at (0,0) { \includegraphics[width=0.3\linewidth]{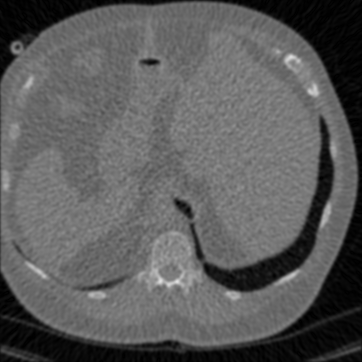}};
\spy on (1.35,.85) in node [left] at (.3,.75);
\node at (-.8,2.5){\footnotesize{PSNR 29.96}};
\node at (-.8,2){\footnotesize{SSIM 0.713}};
\end{tikzpicture}
\hfill
\begin{tikzpicture}[spy using outlines={rectangle,white,magnification=2,size=1.5cm, connect spies}]
\node[anchor=south west,inner sep=0]  at (0,0) { \includegraphics[width=0.3\linewidth]{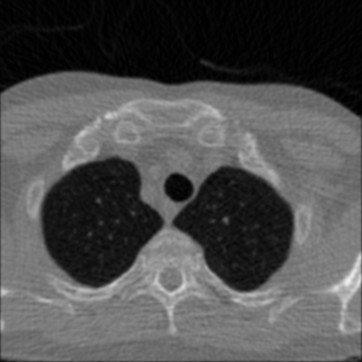}};
\spy on (1.35,.55) in node [left] at (.3,.75);
\node at (-.8,2.5){\footnotesize{PSNR 31.37}};
\node at (-.8,2){\footnotesize{SSIM 0.787}};
\end{tikzpicture}
%\end{subfigure}%
%\vspace{0.2cm}
%\begin{subfigure}[t]{.5\textwidth}
%\caption*{}
\end{subfigure}%

\begin{subfigure}[t]{.5\textwidth}
\centering
  \caption*{\textbf{PD3O}  $\boldsymbol{\lambda}$} 
  \begin{tikzpicture}[spy using outlines={rectangle,white,magnification=2,size=1.5cm, connect spies}]
\node[anchor=south west,inner sep=0]  at (0,0) { \includegraphics[width=0.3\linewidth]{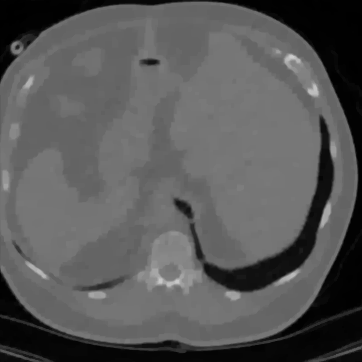}};
\spy on (1.35,.85) in node [left] at (.3,.75);
\node at (-.8,2.5){\footnotesize{PSNR 34.72}};
\node at (-.8,2){\footnotesize{SSIM 0.895}};
\end{tikzpicture}
\hfill
%\end{subfigure}%
%\vspace{0.2cm}
%\begin{subfigure}[t]{.14\textwidth}
%\caption*{}
  \begin{tikzpicture}[spy using outlines={rectangle,white,magnification=2,size=1.5cm, connect spies}]
\node[anchor=south west,inner sep=0]  at (0,0) { \includegraphics[width=0.3\linewidth]{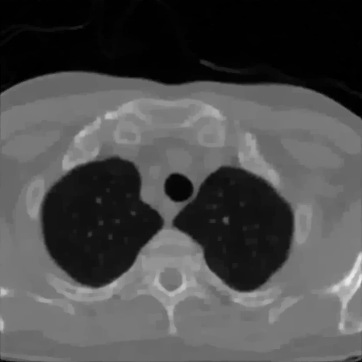}};
\spy on (1.35,.55) in node [left] at (.3,.75);
\node at (-.8,2.5){\footnotesize{PSNR 34.21}};
\node at (-.8,2){\footnotesize{SSIM 0.899}};
\end{tikzpicture}
\end{subfigure}%

\begin{subfigure}[t]{.5\textwidth}
\centering
  \caption*{\textbf{PD3O}  $\boldsymbol{\Lambda}_{\Theta}$} 
  \begin{tikzpicture}[spy using outlines={rectangle,white,magnification=2,size=1.5cm, connect spies}]
\node[anchor=south west,inner sep=0]  at (0,0) { \includegraphics[width=0.3\linewidth]{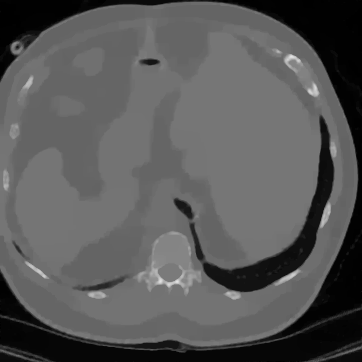}};
\spy on (1.35,.85) in node [left] at (.3,.75);
\node at (-.8,2.5){\footnotesize{PSNR 35.41}};
\node at (-.8,2){\footnotesize{SSIM 0.901}};
\end{tikzpicture}
\hfill
%\end{subfigure}%
%\vspace{0.2cm}
%\begin{subfigure}[t]{.14\textwidth}
%\caption*{}
  \begin{tikzpicture}[spy using outlines={rectangle,white,magnification=2,size=1.5cm, connect spies}]
\node[anchor=south west,inner sep=0]  at (0,0) { \includegraphics[width=0.3\linewidth]{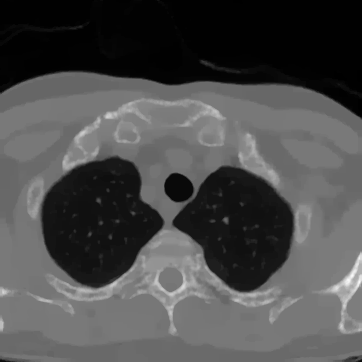}};
\spy on (1.35,.55) in node [left] at (.3,.75);
\node at (-.8,2.5){\footnotesize{PSNR 35.85}};
\node at (-.8,2){\footnotesize{SSIM 0.918}};
\end{tikzpicture}
\end{subfigure}%

\begin{subfigure}[t]{.5\textwidth}
\centering
  \caption*{$\boldsymbol{\Lambda_{\Theta}}$}  
  \begin{tikzpicture}[spy using outlines={rectangle,white,magnification=2,size=1.5cm, connect spies}]
\node[anchor=south west,inner sep=0]  at (0,0) {\includegraphics[width=0.3\linewidth]{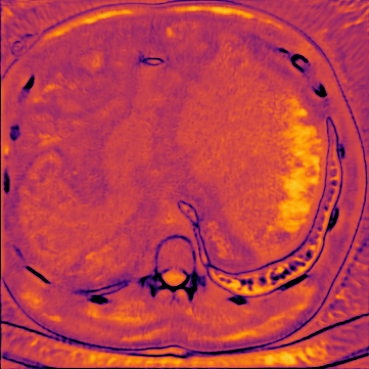}};
\spy on (1.35,.85) in node [left] at (.3,.75);
\node at (-.8,2){\phantom{\footnotesize{SSIM 0.901}}};
\end{tikzpicture}
\hfill
%\end{subfigure}%
%\vspace{0.2cm}
%\begin{subfigure}[t]{.14\textwidth} 
%\caption*{}
  \begin{tikzpicture}[spy using outlines={rectangle,white,magnification=2,size=1.5cm, connect spies}]
\node[anchor=south west,inner sep=0]  at (0,0) {\includegraphics[width=0.3\linewidth]{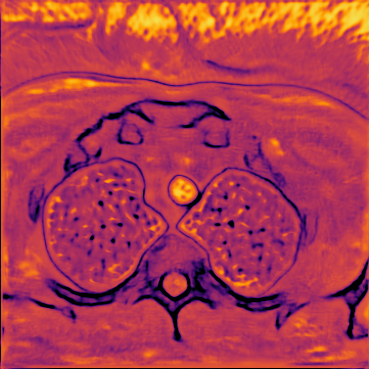}};
\spy on (1.35,.55) in node [left] at (.3,.75);
\node at (-.8,2){\phantom{\footnotesize{SSIM 0.901}}};
\end{tikzpicture}
\end{subfigure}%

\includegraphics[width=\linewidth]{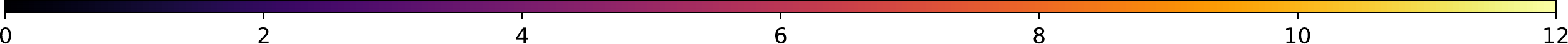}

\vspace{-0.4cm}

\begin{subfigure}[t]{.5\textwidth}
\centering
  \caption*{\textbf{Target}}  
  \begin{tikzpicture}[spy using outlines={rectangle,white,magnification=2,size=1.5cm, connect spies}]
\node[anchor=south west,inner sep=0]  at (0,0) {\includegraphics[width=0.3\linewidth]{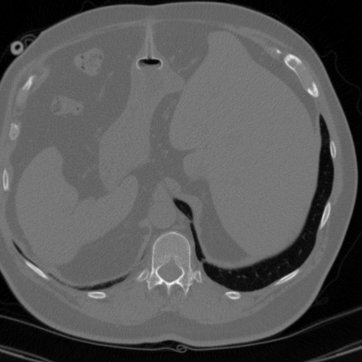}};
\spy on (1.35,.85) in node [left] at (.3,.75);
\node at (-.8,2){\phantom{\footnotesize{SSIM 0.901}}};
\end{tikzpicture}
\hfill
%\end{subfigure}%
%\vspace{0.2cm}
%\begin{subfigure}[t]{.14\textwidth}
%\caption*{}
  \begin{tikzpicture}[spy using outlines={rectangle,white,magnification=2,size=1.5cm, connect spies}]
\node[anchor=south west,inner sep=0]  at (0,0) {\includegraphics[width=0.3\linewidth]{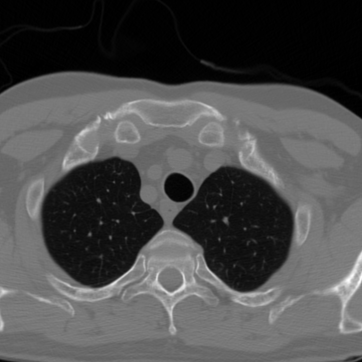}};
\spy on (1.35,.55) in node [left] at (.3,.75);
\node at (-.8,2){\phantom{\footnotesize{SSIM 0.901}}};
\end{tikzpicture}
\end{subfigure}%

\caption{Different reconstructions obtained with PD3O employing a scalar regularization and the regularization parameter-maps obtained with the proposed CNN. From top to bottom: initial FBP-reconstruction, PD3O $\lambda$, PD3O $\Lambda_{\Theta}$, spatial regularization parameter-map and ground truth image. The corresponding PSNR and SSIM values are given in the top left corner of the image.} \label{fig:CT_img_lambdamaps}
\end{figure}
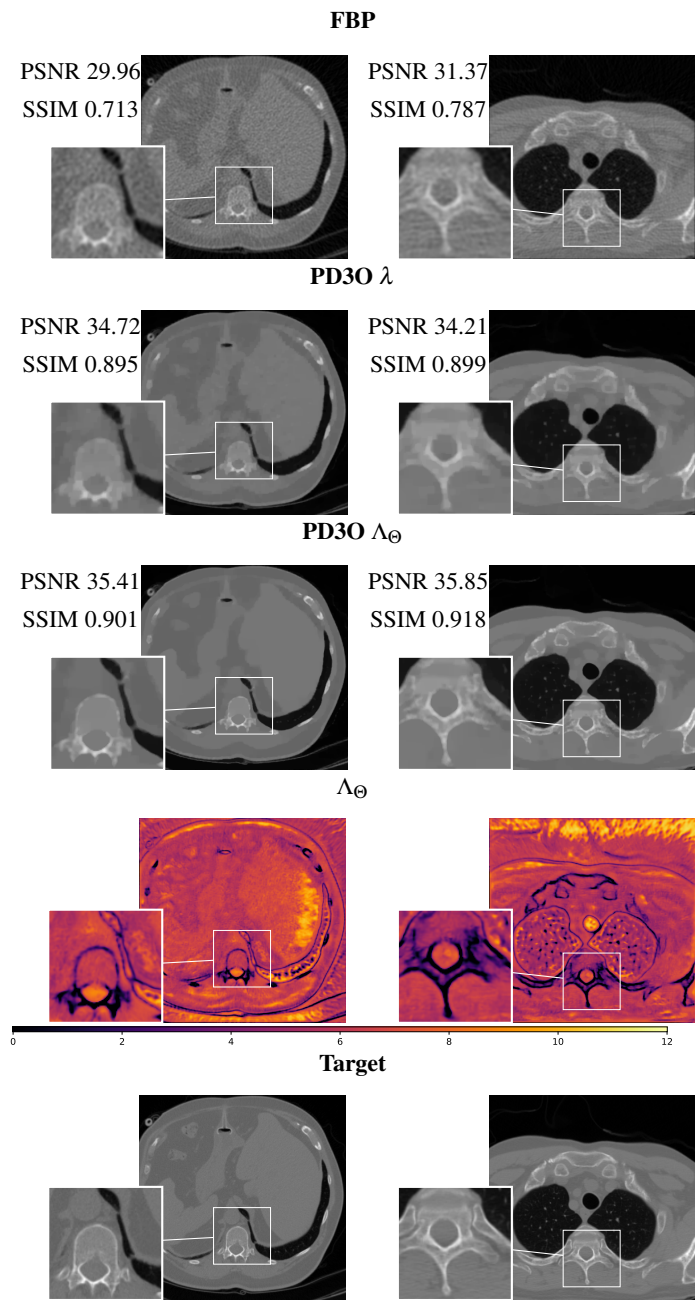

In Figure~\ref{fig:CT_img_lambdamaps} we compare the FBP reconstruction with the PD3O reconstructions where we use (i) a scalar parameter ($\lambda$), chosen to maximize the PSNR ``on average'', and (ii) our computed spatially dependent parameter map ($\boldsymbol\Lambda_{\Theta}$). Using the latter, we obtain a significant improvement in the reconstruction both visually and in terms of quality measures, e.g., PSNR and SSIM.  In particular, sharp edges are retained, while the constant regularizing parameter
results in blurry and blocky-like reconstructions. One can observe that the network attributes higher regularization parameters to image content with smooth structures while it yields lower regularization parameters at the edges to prevent smoothing.

\section{Discussion}

We have presented a data-driven approach to automatically estimate  spatially dependent parameter-maps for TV regularization for the low dose X-ray CT tomography reconstruction problem. Although only the TV regularization is considered in this paper, higher order  or combinations of regularizers can be used for different CT applications, see \cite{Papoutsellis_2021, Warr2021}. Moreover, our unrolled  framework
is quite flexible and can be easily used for other modalities such as qualitative and quantitative MRI reconstruction, image denoising  as well as their dynamic versions, see \cite{kofler2023learning}. Finally,  more sophisticated network architectures than the U-Net have been proposed recently, e.g. \cite{lu2022transformer, liang2021swinir}, which could be potentially adopted for the estimation of the regularization parameter-maps as well.

 \section{Acknowledgments}
 The authors acknowledge the support of the German Research Foundation (DFG) under Germany's Excellence Strategy – The Berlin Mathematics
 Research Center MATH+ (EXC-2046/1, project ID: 390685689) as this work was initiated during the Hackathon event “Maths Meets Image”, Berlin, March 2022, which was part of the MATH+ Thematic Einstein Semester on ``Mathematics of Imaging in Real-World Challenges". This work is further supported  via the MATH+ project EF3-7.

 This work was funded by the UK EPSRC grants the ``Computational Collaborative Project in Synergistic Reconstruction for Biomedical Imaging" (CCP SyneRBI) EP/T026693/1; ``A Reconstruction Toolkit for Multichannel CT" EP/P02226X/1 and ``Collaborative Computational Project in tomographic imaging’' (CCPi) EP/M022498/1 and EP/T026677/1. This work made use of computational support by CoSeC, the Computational Science Centre for Research Communities, through CCP SyneRBI and CCPi.

 This work is part of the Metrology for Artificial Intelligence for Medicine (M4AIM) project that is funded by the German Federal Ministry for Economic Affairs and Energy (BMWi) in the framework of the QI-Digital initiative.

%-------------------------------------------------------------------------------------------
% \printbibliography
% \bibliography{refs.bib}

% \end{document}

%-------------------------------------------------------------------------------------------
\printbibliography

\end{document}